\documentclass[12pt]{amsart}

\usepackage{epsfig}
\usepackage{verbatim}

\begin{document}
\title {On Seymour's Decomposition Theorem}
\maketitle

\begin {center}
S. R. Kingan 
\footnote{The author is partially supported by  PSC-CUNY grant number 66305-00 43.} \\     
Department of Mathematics \\
Brooklyn College, City University of New York\\
 Brooklyn, NY 11210\\
skingan@brooklyn.cuny.edu\\  
\end {center}

\begin{abstract}  Let $\mathcal M$ be a class of matroids closed under minors and isomorphism. Let $N$ be a matroid in $\mathcal M$ with an exact $k$-separation $(A, B)$. We say $N$ is a $k$-decomposer for $\mathcal M$ having $(A, B)$ as an inducer, if every matroid $M\in \mathcal M$ having $N$ as a minor has a $k$-separation $(X, Y)$ such that, $A\subseteq X$ and $B\subseteq Y$. Seymour [3, 9.1] proved that a matroid $N$ is a $k$-decomposer for an excluded-minor class, if certain conditions are met for all 3-connected matroids $M$ in the class, where $|E(M)-E(N)|\le 2$. We reinterpret Seymour's Theorem in terms of the connectivity function and give a check-list that is easier to implement because case-checking is reduced.  

\end{abstract}

\bigskip 

\section {Introduction}

Excluded minor results are obtained by excluding certain minors and determining the 3-connected matroids without the excluded minors. If the class has few infinite families, then it may be possible to identify all of them.   Often, however, there are too many 3-connected matroids with no  pattern to get a precise identification of them. One approach is to raise the connectivity. It turns out that while 3-connected matroids are plentiful, internally 4-connected matroids are fewer. Determining all the internally 4-connected matroids in a class is not the same as determining all the 3-connected matroids in the class, but it gives a good understanding of the structure of the class. This amounts to asking whether a non-minimal exact 3-separation in a 3-connected matroid $N$ gets carried forward in all 3-connected matroids in the class that contain $N$ as a minor. Let us call such a matroid a $3$-decomposer because it captures, in a sense, the decomposition present in all 3-connected matroids containing it as a minor. What is needed then is a quick and easy check to determine if a matroid is a 3-decomposer. In this paper we reinterpret on Seymour's Decomposition Theorem [3, 9.1] to present such an easy check.

The matroid terminology follows Oxley [2]. If $M$ and $N$ are matroids on the sets $E\cup x$ and $E$ where $x   \not\in E$, then $M$ is a  {\it single-element extension} of $N$ if $M \backslash x = N$,  and $M$ is a {\it single-element coextension} of $N$ if $M/x=N$. If $N$ is a 3-connected matroid, then an extension $M$ of $N$ is  3-connected provided $x$ is not in a 1- or 2-element circuit of $N$. Likewise, $M$ is a 3-connected coextension of $N$ provided $x$ is not in a 1- or 2-element cocircuit of $N$.

Let $M$ be a matroid and $X$ be a subset of the ground set $E$. The {\it connectivity function}  $\lambda$ is defined as $\lambda (X) = r(X) + r(E-X) - r(M)$.  Observe that $\lambda (X) = \lambda (E-X)$.  For $k\ge 1$, a partition $(A, B)$ of $E$ is called a $k$-separation if $|A|\ge k$, $|B|\ge k$, and $\lambda (A) \le k-1$.  When $\lambda (A)=k-1$, we call $(A, B)$ an {\it exact k-separation}.  When $\lambda (A)=k-1$ and $|A|=k$ or $|B|=k$ we call $(A, B)$ a {\it minimal exact k-separation}. For $n\ge 2$, we say $M$ is {\it n-connected} if $M$ has no $k$-separation for $k\le n-1$.  A  $k$-connected matroid is {\it internally $(k+1)$-connected} if it has no non-minimal exact $k$-separations. 
In particular, a simple matroid is 3-connected if $\lambda (A)\ge 2$ for all partitions $(A, B)$ with $|A|\ge 3$ and $|B|\ge 3$. A 3-connected matroid is {\it internally $4$-connected}  if $\lambda (A)\ge 3$ for all partitions $(A, B)$ with $|A|\ge 4$ and $|B|\ge 4$. In this case $\lambda (A)= 2$ is allowed only when either  $|A|$ or $|B|$ has size at most 3.

Let $\mathcal M$ be a class of matroids closed under minors and isomorphisms. Let $k\ge 1$ and $N$ be a matroid in $\mathcal M$ having an exact  $k$-separation $(A, B)$. Let $M$ be a matroid in  $\mathcal M$ having an $N$-minor.  
We say $N$ is a {\it $k$-decomposer} for $M$ having $(A, B)$ as an {\it inducer} if $M$ has an exact $k$-separation $(X, Y)$ such that $A\subseteq X$ and $B\subseteq Y$. From a practical point of view we are interested in 3-connected matroids with non-minimal exact 3-separations. 
Suppose $M$ is a 3-connected matroid having a 3-connected minor $N$ and $N$ has a non-minimal exact 3-separation $(A, B)$. If $N$ is a 3-decomposer for $M$, then $M$ has a $3$-separation $(X, Y)$ such that $A\subseteq X$ and $B\subseteq Y$. In this case $\lambda (X)=2$ and $|X|\ge 4$ and $|Y|\ge 4$.  Hence, if $N$ is a 3-decomposer for $M$, then $M$ is not internally 4-connected. The converse is not true. 

The first excluded minor class of matroids characterized in this manner was the class of regular matroids [3]. A matroid is {\it regular} if it has no minor isomorphic to the Fano matroid $F_7$ or its dual $F^*_7$. To decompose regular matroids, Seymour developed the  Splitter Theorem, the Decomposition Theorem, and the notion of 3-sums. The Splitter Theorem describes how 3-connected matroids can be systematically built-up and the Decomposition Theorem (the subject of this paper) describes the conditions under which an exact 3-separation in a matroid gets carried forward to all matroids containing it. 

The next theorem is the main result in this paper.  It is equivalent to Seymour's Decomposition Theorem [3, Theorem 9.1]. When Seymour developed his Decomposition Theorem in 1980 finding the connectivity function of a matroid was a daunting task. Now the connectivity function of a matroid can be obtained with a click of a button in any matroid software. As such it is only natural to reinterpret Seymour's conditions in terms of the connectivity function to make it easier to use. Our theorem is, in a sense, a back-to-basics result because it says Seymour's Theorem is still a good way of finding internally 4-connected matroids.
\bigskip

\noindent{\bf Theorem 1.1.} {\it  Let $N$ be a simple and cosimple matroid in $\mathcal{M}$ with an exact $k$-separation $(A,B)$, such that $A$ is the union of circuits and the union of cocircuits.  Suppose $M\in \mathcal {M}$.

\begin{enumerate}
\item[(i)] If $M$ is a simple single-element extension of $N$ such that $M\backslash e=N$, then $\lambda_M(A)=k-1$ or $\lambda_M(A\cup e)=k-1$.

\item[(ii)] If $M$ is a cosimple single-element coextension of $N$ such that $M/ f=N$, then $\lambda_M(A)=k-1$ or $\lambda_M(A\cup f)=k-1$.

\item[(iii)] Each matroid $M$ that is a cosimple single-element coextension of a Type (i) matroid or a simple single-element extension of a Type (ii) matroid satisfies one of the following conditions:

\begin{enumerate}
\item[(a)] $\lambda _{M/f}(A)=k-1$ and $\lambda _{M\backslash e}(A)=k-1$;
\item[(b)] If $\lambda _{M/f}(A)=k-1$ and $\lambda _{M\backslash e}(A\cup f)=k-1$, then either $\lambda _M(A\cup f)=k-1$ or $\{e, f, g\}$ is a triad or triangle with $g\in A$;
\item[(c)]If $\lambda _{M/f}(A\cup e)=k-1$ and $\lambda _{M\backslash e}(A)=k-1$, then either $\lambda _M(A\cup e)=k-1$ or $\{e, f, g\}$ is a triad or triangle with $g\in A$;  or 
\item[(d)] If $\lambda _{M/f}(A\cup e)=k-1$ and $\lambda _{M\backslash e}(A\cup f)=k-1$, then $\{e, f, g\}$ is a triangle or triad in $M$ with $g\in A$.
\end{enumerate}
\end{enumerate}
\noindent Then $N$ is a $k$-decomposer for every matroid in $\mathcal{M}$ with an $N$-minor.}

\bigskip

There are several points to note. First, observe that for $k\ge 4$, we do not have a splitter-type theorem that tells us how to construct $k$-connected matroids. From a practical point of view this result is used only for 3-connected matroids and 3-decomposers. 
Second, observe that if $N$ is 3-connected, a simple single-element extension and a cosimple single-element coextension are also 3-connected. 
 
Third, computing a matroid $M$ where $M\backslash e/f=N$ is clearly more work than computing single-element extensions and coextensions of $N$, especially since we have to consider every single column and row and cannot use isomorphism.  Anything that can be done to reduce the computation is worth it. Due to (iii)(a), many of the choices for extension columns and coextension rows are eliminated. This is a considerable savings of computation. 
 
Lastly, observe that there is no symmetry between $A$ and $B$. So for example, in the case of a matroid $M$ of the form $M\backslash e/f=N$, having $\lambda_{M}(A)=2$ is permissable, but $\lambda (A\cup \{e, f\})=\lambda _M(B)=2$ could be a problem. An example is given to illustrate this point at the end of the paper. 

This example, also serves to show that if $N$ has a non-minimal exact two separation $(A, B)$ and for all $M$, such that $|E(M)-E(N)|\le 2$, $M$ has a non-minimal exact 3-separation $(X, Y)$ such that $A\subseteq X$ and $B\subseteq Y$, that does not mean $N$ is a 3-decomposer. In other words, just a 2-element check will not do, not even for binary matroids. The highly specialized conditions in Seymour's theorem (and equivalently our interpretation) must hold.  

Next, we examine the situation when we add additional hypothesis on the exact 3-separation $(A, B)$ of $N$. Specifically, suppose $A$ is not only a union of circuits and a union of cocircuits, but a circuit and a cocircuit. Then  $r_N(A)=|A|-1$ and since $B$ is a hyperplane $r_N(B)=r(N)-1$
\begin{eqnarray*}
 \lambda _N\left(A\right) & = & r_N  \left(A\right)  +r_N\left(B\right) - r \left(M\right)    \\
                  & = & |A|-1+r\left(M\right)-1- r\left(N\right)\\
                  & = & |A|-2 
\end{eqnarray*} 
 
Now, suppose further that $A$ is a 4-element circuit and cocircuit. Then $\lambda _N (A)=2$.  This is precisely the condition that is in  Mayhew, Royle, and Whittle's Decomposition Theorem [1, Lemma 2.10].   

\bigskip

\noindent {\bf Corollary 1.2.} {\it (Mayhew, Royle, and Whittle) Let $N$ be a $3$-connected matroid in $\mathcal{M}$ such that $N$ is not a wheel or a whirl and $N$ has a non-minimal exact $3$-separation $(A,B)$, where $A$ is a $4$-element circuit and cocircuit. If $A$ is a circuit and a cocircuit in every $3$-connected single-element extension and  coextension of $N$ in $\mathcal M$, then $N$ is a $3$-decomposer for every $3$-connected matroid in $\mathcal{M}$ having an $N$-minor.}
\bigskip

We prove Theorem 1.1 by showing that it is equivalent to Seymour's Decomposition Theorem [3, 9.1]. To do this, in Section 2  we make a slight modification of Seymour's Theorem and show that the modified version is equivalent to the original version in [3]. In Section 3 we describe the updated techniques for finding single-element extensions and the connectivity function that led to this perspective on decomposition. In Section 4 we  give the proof of Theorem 1.1 and Corollary 1.2. Finally, in Section 5 we illustrate Theorem 1.1. by giving a short explanation for why $R_{12}$ is a 3-decomposer for regular matroids. 


\bigskip  

\section {A modification of Seymour's Decomposition Theorem}

In this section we make a slight modification of  Seymour's Decomposition Theorem [3, Theorem 9.1].   

\bigskip

\noindent {\bf Theorem 2.1.} {\it (Seymour's Decomposition Theorem) Let $\mathcal M$ be a class of matroids, closed under minors and isomorphism. Let $N\in \mathcal M$, and let $(A, B)$ be a partition of $E(N)$, with $\lambda _N (A, B)=k$. Suppose $N, \mathcal M$ have the following properties:
\begin{enumerate}
\item [(i)] For each $x\in A$, there is a circuit $C$ and a cocircuit $D$ of $N$ containing $x$, with $C, D \subseteq A$.

\item [(ii)] For each $M\in \mathcal M$, if $M\backslash e = N$ and $e$ is not a coloop of $M$, there is a circuit $C_e$ of $M$ with $e\in C_e$ and $C_e - \{e\}$ included in one of $A$, $B$.

\item [(iii)] For each $M\in \mathcal M$, if $M/ f = N$ and $f$ is not a loop of $M$, there is a cocircuit $D_f$ of $M$ with $f\in D_f$ and $D_f - \{f\}$ included in one of $A$, $B$.

\item [(iv)] For each $M\in \mathcal M$, if $M\backslash e/ f = N$, suppose there is a cocircuit $D$ of $M$ with $\{e, f\}\subset D\subseteq B\cup \{e, f\}$; then either there is a circuit $C$ of $M$ with $x\in C\subseteq B\cup \{e, f\}$ or $e$ is parallel to an element of $A$ in $M/f$. 

\item [(v)] For each $M\in \mathcal M$, if $M\backslash e/ f = N$, suppose there is a circuit $C$ of $M$ with $\{e, f\}\subset C\subseteq B\cup \{e, f\}$; then either there is a cocircuit $D$ of $M$ with $f\in D\subseteq B\cup \{e, f\}$ or $f$ is in series with an element of $A$ in $M\backslash e$. 

\end{enumerate}

Then $\lambda _M (A, B)=k$ for each $M\in \mathcal M$ with $N$ as a minor.
}
\bigskip

The next result is equivalent to Theorem 2.1, but emphasizes the role of simple single-element extensions and cosimple single-element coextensions. As noted earlier, if $N$ is 3-connected, a simple single-element extension and a cosimple single-element coextension are also 3-connected.  

\bigskip

\noindent{\bf Theorem 2.2.} {\it Let $N$ be a simple and cosimple matroid in $\mathcal{M}$ with an exact $k$-separation $(A,B)$ such that $A$ is the union of circuits and the union of cocircuits.  Suppose $M\in \mathcal {M}$.
\begin{enumerate}
\item[(i)] If $M$ is a simple single-element extension of $N$ by element $e$, then $M$ has a circuit $C_e$ containing $e$ such that $C_e$ is contained in $A\cup e$ or $B\cup e$.

\item[(ii)] If $M$ is a cosimple single-element extension of $N$ by element $f$, then $M$ has a cocircuit $D_f$ containing $f$ such that $D_f$ is contained in $A\cup f$ or $B\cup f$.

\item[(iii)] For each matroid $M$ that is a cosimple single-element coextension of matroids of Type (i) or a simple single-element extension of matroids of Type (ii) having a circuit or cocircuit $R$ containing $e, f$ and $R\subseteq B\cup \{e, f\}$

\begin{enumerate}
\item[(a)] If $R$ is a circuit, then either $M$ has a cocircuit $D_f$ such that $f\in D_f\subseteq B\cup \{e, f\}$ or there is an element $g\in A$ such that $\{e, f, g\}$ is a triad in $M$. 
\item[(b)] If $R$ is a cocircuit, then either $M$ has a circuit $C_e$ such that $e\in C_e \subseteq B\cup \{e, f\}$ or there is an element $g\in A$ such that $\{e, f, g\}$ is a triangle in $M$. 
\end{enumerate}

\end{enumerate}
\noindent Then $N$ is a $k$-decomposer for every matroid in $\mathcal{M}$ with an $N$-minor.}

\bigskip

\noindent {\bf Proof.} The result follows provided we verify conditions (i) to (v) of Theorem 3.1. Note that (i) follows from the hypothesis because  $A$ is the union of circuits and the union of cocircuits. In order to verify Theorem 3.1(ii) and (iii), by duality we need to establish only Theorem 3.1(ii). Suppose $M\backslash e=N$ and $e$ is not a coloop of $M$. 

If $M$ is not simple, then $C=C'\cup e$ is a circuit of $M$, for some $C'\subseteq E(N)$ such that $|C'|\le 1$. As $C'\subseteq A$ or $C'\subseteq B$, it follows that $C'-e$ is included in $A$ or $B$. Thus we have Theorem 3.1(ii). 

If $M$ is simple, then since $M$ is a simple single-element extension of $N$ by element $e$, $(A\cup e,B)$ is an exact $k$-separation for $M$. Therefore $r_M(A)=r_M(A\cup e)$. So there is a circuit $C_e$ of $M$ with $e\in C_e$ and $C_e-\{e\}\subseteq A$. Thus we have Theorem 3.1(ii). 

In order to show Theorem 3.1(iv) and (v) again by duality we need to establish only (iv). Suppose  $M\backslash e/f=N$ and $M$ has a circuit $R$ such that $\{e,f\}\subset R\subseteq B\cup\{e,f\}$. 
First, consider the case that $M$ is not simple nor cosimple. Then because $N$ is simple and cosimple, $e$ or $f$ is contained in a non-trivial series or parallel class $T$ of $M$. 

Suppose $T$ is a series class of $M$. If  $f\in T$, then by orthogonality with $R$, there is $g\in (R-f)\cap T$ and therefore $D_f=\{f,g\}$ is a cocircuit of $M$ containing $f$ and $D_f\subseteq E\cup \{e, f\}-X$. This gives us Theorem 3.1(iv). If $f\not \in T$, then $e\in T$. In this case $T$ is a series class of $M/f$. So $T-e$ is a non-empty set of coloops of $N=M/f\backslash e$; a contradiction to the hypothesis.

Suppose $T$ is a parallel class of $M$.
If $e\not\in T$, then $f\in T$ and $T-f$ is a non-empty set of loops of $N=M/f\backslash e$; a contradiction. Therefore $e\in T$. By (ii) applied to $M\backslash e$, there is a cocircuit $D_f$ of $M\backslash e$ containing $f$ such that $D_f-\{f\}$ is contained in one of $A$ or $B$. Let $D'_f$ be a cocircuit of $M$ such that $D'_f-e=D_f$. If $D_f-f\subseteq B$, then $f$ is contained in $D'_f$ and $D'\subseteq B\cup\{e,f\}$. Theorem 3.1(iv) follows.

 Therefore we may assume $D'_f-f\subseteq A$. By orthogonality with $R$, $D'_f=D_f\cup e$, since $D_f\cap R=\{f\}$. If $g\in T-e$, then, by orthogonality with $D'_f$, $g\in A$ or $g=f$. We will prove that $T=\{e,f\}$. Suppose $T \neq \{e,f\}$. We can choose $g$ such that $g\in A$ and  $(R-f)\bigtriangleup\{e,g\}$ contains a circuit $R'$ of $N=M/f\backslash e$ satisfying $g\in R'$. So $R'-B=\{g\}$; a contradiction to orthogonality because $N$ has a cocircuit containing $g$ and contained in $A$ (recall that $(A,B)$ is a special $k$-separation for $N$). Thus we have shown that $T=\{e,f\}$.  This is a contradiction because $T$ is a parallel class and therefore $T=R$. From now on we may assume $M$ is both simple  and cosimple. 

Suppose $M$ is a cosimple single-element coextension of a Type (i) matroid or a simple single-element extension of a type (ii) matroid. Suppose this is not the case. Then $M/f$ is not simple and $M\backslash e$ is not cosimple because $M$ is both simple and cosimple. So $M$ has a triangle $T$ and a triad $T^*$ such that $f\in T$ and $e\in T^*$. As $N=M\backslash e/f$ is both simple and cosimple, it follows that $\{e,f\}\subseteq T\cap T^*$, say $T^*=\{e,f,g\}$. If $g\in A$, then $e$ and $f$ are in series in $M\backslash e$ and so Theorem 3.1(iv) follows. Thus we may assume that $g\in B$. Again Theorem 3.1(iv) follows because $f\in T^*\subseteq B\cup\{e,f\}$. Hence proved. $\qed$
\bigskip


\section {Updated Techniques}

The author was introduced to structural results that involved computing single-element extensions by James Oxley. In the past, determining the non-isomorphic single-element extensions of a $GF(q)$-representable matroid required geometric insights and lengthy arguments. Now the author's matroid software program, Oid, calculates single-element extensions and groups them into isomorphism classes at the click of a button. 

Let $N$ be a $GF(q)$-representable $n$-element rank-$r$  matroid  represented by the matrix $A=[I_r|D]$ over $GF(q)$. The columns of A may be viewed as a subset of the columns of the matrix that represents the projective geometry $PG(r - 1, q)$. Let $M$ be a simple single-element extension of $N$ over $GF(q)$. Then $N=M\backslash e$ and $M$ may be represented by $[I_r|D']$, where $D'$ is the same as $D$, but with one additional column corresponding to the element $e$. The new column is distinct from the existing columns and has at least two non-zero elements. If the existing columns are labeled $\{1, \dots , r, \dots , n\}$, then the new column is labeled $(n+1)$.

Suppose $M$ is a cosimple single-element coextension of $N$ over $GF(q)$. Then $N=M/f$ and $M$ may be represented by the matrix $[I_{r+1}| D'']$, where $D''$ is the same as $D$, but with one additional row. The new row is distinct from the existing rows and has at least two non-zero elements. Since the columns of the original matrix $A=[I_r|D]$ are labeled $\{1, \dots r, \dots , n\}$, the columns of $[I_{r+1}| D'']$ are labeled 
$\{1, \dots , r+1, \dots , n+1\}$. The coextension element $f$ corresponds to column $r+1$. The coextension row is selected from $PG(n-r, q)$, which means there could be a much larger selection of row vectors for the coextension. 

We refer to the simple single-element extensions of $N$ as Type (i) matroids and the cosimple single-element coextensions of $N$ as Type (ii) matroids. The structure of type (i) and Type (ii) matroids are shown in Figure 1.

\begin{figure}[h]
\centering
\epsfxsize 5in \epsfbox{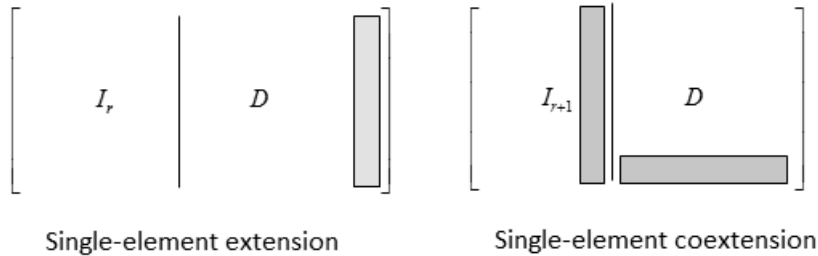}
\caption{Structure of Type (i) and Type (ii) matroids }
\end{figure}

We can visualize the new element $f$  as appearing in the new dimension and lifting several points into the higher dimension. Observe that $f$ forms a cocircuit with the elements corresponding to the non-zero elements in the new row.  Note that in $[I_{r+1}|D'']$ the labels of columns beyond $r$ are increased by 1 to accomodate the new column $r+1$. So, if in a rank 6 matroid $N$, $X=\{1, 2, 5, 7, 10, 11\}$, then in its single-element coextension $M$, the set corresponding to $X$ is $X'=\{1, 2, 5, 8, 11, 12\}$. 

For example, binary matrix representations for $R_{10}$ and $PG(4, 2)$ are given by the matrices $A$ and $P$ shown below.  

\tiny

\[
A=\left[ 
\begin{array}{ccccc|ccccc}
&&&&&    1&0&0&1&1 \\
&&&&&    1&1&0&0&1 \\
&&I_5&&& 1&1&1&0&0 \\
&&&&&    0&1&1&1&0 \\
&&&&&    0&0&1&1&1
\end{array} 
\right] 
\]

\[
P=\left[ 
\begin{array}{ccc|ccccccccccccccccccccccccccccccc}
&&&    1&0&0&0&0&0&0&0&0&0&0&0&0&0&0&0&1&1&1&1&1&1&1&1&1&1&1&1&1&1&1 \\ 
&&&    0&1&0&0&0&0&0&0&0&1&1&1&1&1&1&1&0&0&0&0&0&0&0&1&1&1&1&1&1&1&1 \\ 
&I_5&& 0&0&1&0&0&0&1&1&1&0&0&0&1&1&1&1&0&0&0&1&1&1&1&0&0&0&0&1&1&1&1  \\ 
&&&    0&0&0&1&0&1&0&1&1&0&1&1&0&0&1&1&0&1&1&0&0&1&1&0&0&1&1&0&0&1&1 \\ 
&&&    0&0&0&0&1&1&1&0&1&1&0&1&0&1&0&1&1&0&1&0&1&0&1&0&1&0&1&0&1&0&1\\ 
\end{array} 
\right] 
\] 
\normalsize

\noindent Comparing $A$ and $P$ and adding to $A$ the columns in $P$ missing in 
$A$ gives us two isomorphism classes. Each class contains several choices for the new column representing the element $e$ in the single-element extension of $R_{10}$. So, for instance, in Figure 2, the last columm $[0 0 0 1 1]$ of Extension 1 may be replaced by any one of the columns in the list to get an isomorphic matroid.

\begin{figure}[h]
\centering
\epsfxsize 4in \epsfbox{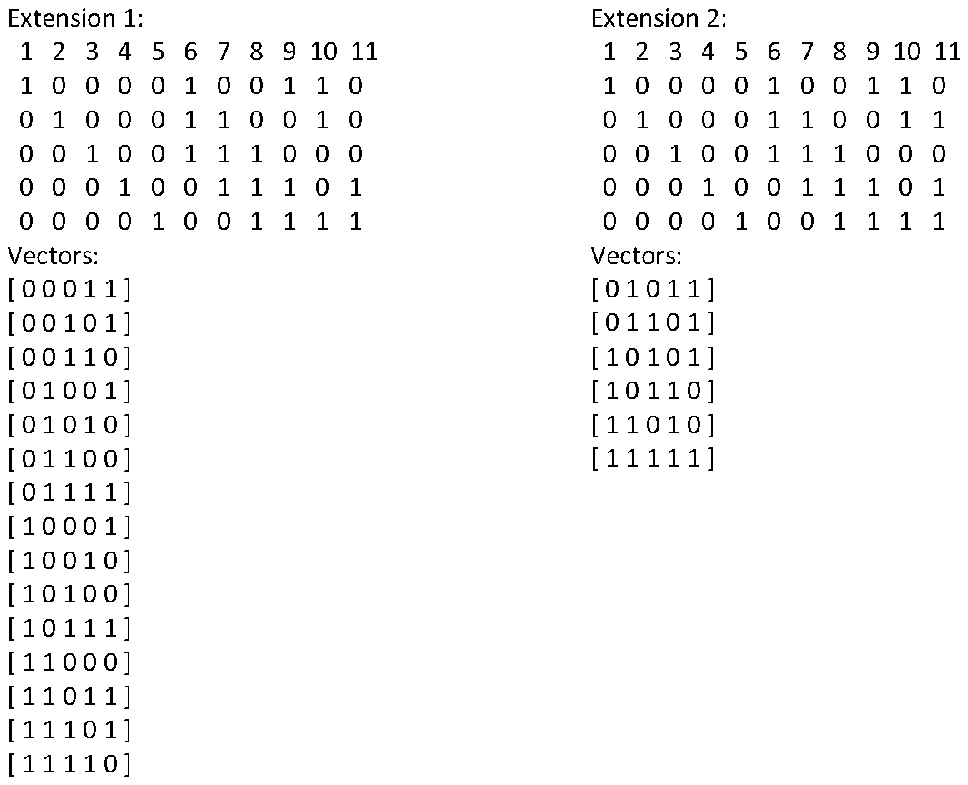}
\caption{The binary simple single-element extensions of $R_{10}$ }
\end{figure}

A binary matrix representation for $R_{12}$ is shown below. A coextension of $R_{12}$ has rank 7 and 13 elements. So the selection of rows to add at the bottom will come from $PG(6, 2)$, which is too large to present here.

\[
B=\left[ 
\begin{array}{cccccc|cccccc}
&&&&&&   1&1&1&0&0&0 \\
&&&&&&   1&1&0&1&0&0 \\
&&&I_6&&&1&0&0&0&1&0 \\
&&&&&&   0&1&0&0&0&1 \\
&&&&&&   0&0&1&0&1&1 \\
&&&&&&   0&0&0&1&1&1 
\end{array} 
\right] 
\]

The single-element coextensions of $R_{12}$ are given in Figure 3. The last row $[0 0 0 1 1]$ of Extension 1 may be replaced by any one of the rows in the list to get an isomorphic matroid.

\begin{figure}[h]
\centering
\epsfxsize 5in \epsfbox{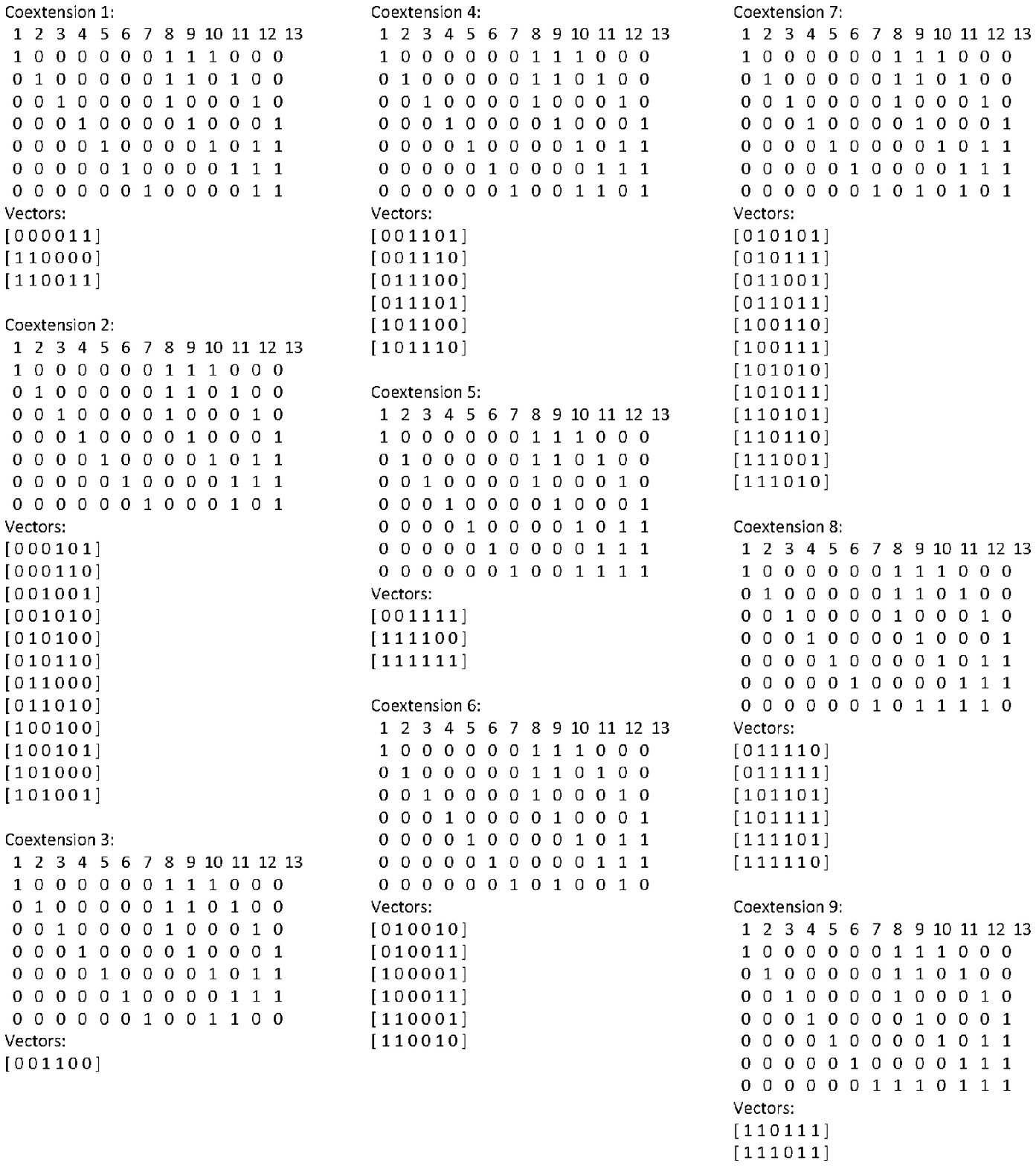}
\caption{The binary cosimple single-element coextensions of $R_{12}$ }
\end{figure}

It is easy to determine if a matroid is non-regular, because a non-regular matroid has an $F_7$ or  $F_7^*$-minor. Only, the first and third coextension of $R_{12}$ are regular. 
The connectivity function is also easily determined. When the connectivity function is known, non-minimal exact 3-separations can be easily found. For example, $\lambda (A)=2$  in $R_{12}$, where $A=\{1, 2, 5, 6, 9, 10\}$.

Once the simple single-element extensions (Type (i) matroids) and cosimple single-element coextensions (Type (ii) matroids) are determined, the number of permissable rows and columns give a bound on the choices for the cosimple single-element extensions of the Type (i) matroids and the simple single-element extensions of the Type (ii) matroids, respectively. 

For example, in $R_{12}$ since only the first and third coextension are regular, the only choices for the cosimple single-element coextensions of Type (i) matroids are the four permissable rows in a Type (ii) matroid with a zero and a one, as well as, the rows that are in series with existing rows with the last entry reversed. Similarly, the only choices for the simple single-element extensions of Type (ii) matroids are the four columns in a Type (i) matroid, as well as, the columns that are in parallel with existing columns with the last entry reversed. The structure of the cosimple single-element coextensions of a Type (i) matroid and the simple single-element extensions of a Type (ii) matroid are shown in Figure 4.
 
\begin{figure}[h]
\centering
\epsfxsize 5in \epsfbox{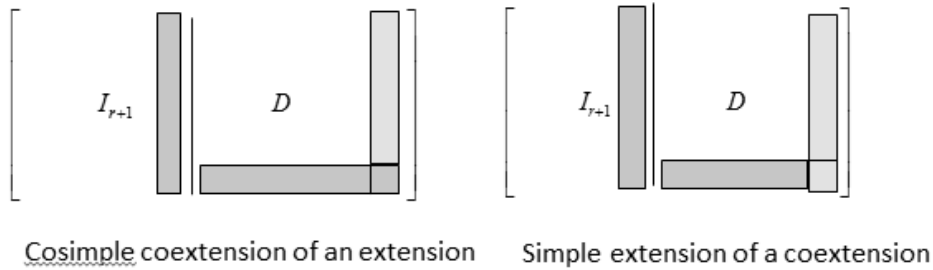}
\caption{Structure of $M$, where $|E(M)-E(N)|=2$ }
\end{figure}


\section {Proof of the Main Theorem}
\bigskip

We begin by proving three lemmas on the connectivity function that are used in the proof  of Theorem 1.1.

\bigskip
\noindent{\bf Lemma 4.1.} {\it Let $N$ be a simple and cosimple matroid in $\mathcal{M}$ with an exact $k$-separation $(A, B)$. Suppose $M\in \mathcal{M}$ such that $M\backslash e = N$ and $e$ is not a loop. Then
\begin{enumerate}
\item [(i)] $\lambda_M(A)=k-1$ if and only if $M$ has a circuit $C_e$ such that $e\in C_e\subseteq B\cup e$;
\item [(ii)] $\lambda_M(A\cup e)=k-1$ if and only if $M$ has a circuit $C_e$ such that $e\in C_e\subseteq A\cup e$.
\end {enumerate}}

\bigskip

\noindent {\bf Proof.}  Since $(A, B)$ is an exact $k$-separation of $N$, $\lambda_M(A)=k-1=\lambda_N(A)$. So, $$r_M(A)+r_M(B\cup e)-r(M)=r_N(A)+r_N(B)-r(N).$$ 
Observe that $r(M)=r(N)$ and $r_M(A)=r_N(A)$. Therefore,  
$$r_M(B\cup e)=r_N(B)...............(1)$$ 
Since $e$ is not a loop, (1) occurs if and only if there is a circuit $C_e$ containing $e$ such that $C_e\subseteq B\cup e$.  This completes the proof of Part (i). The proof of (ii) may be obtained by replacing $A$ with $B$, and noting that $\lambda _M(A\cup e)=\lambda_M(B)=k-1$. $\qed$

\bigskip
\noindent{\bf Lemma 4.2.} {\it Let $N$ be a simple and cosimple matroid in $\mathcal{M}$ with an exact $k$-separation $(A, B)$. Suppose $M\in \mathcal{M}$ such that  $M / f=N$ and $f$ is not a coloop.  Then 
\begin {enumerate}
\item [(i)] $\lambda_M(A)=k-1$ if and only if $M$ has a cocircuit $D_f$ such that $f\in D_f\subseteq B\cup f$;  
\item [(ii)] $\lambda_M(A\cup f)=k-1$ if and only if $M$ has a cocircuit $D_f$ such that $f\in D_f\subseteq A\cup f$. 
\end {enumerate}}
 
\bigskip

\noindent {\bf Proof.}  Again, since  $\lambda_M(A)=k-1=\lambda_N(A)$,
  $$r_M(A)+r_M(B\cup f)-r(M)=r_N(A)+r_N(B)-r(N).$$ 
Observe that $r(M)=r(N)+1$ and $r_M(B\cup f)=r_N(B)+1$. Therefore, 
$$r_M(A)=r_N(A)...............(2)$$ 
But $r_M(A\cup f)=r_N(A)+1$. So $$r_M(A\cup f)=r_M(A)+1.$$
This occurs if and only if there is a hyperplane containing $A$ and avoiding $f$ (since $f$ is not a coloop). It follows that there is a cocircuit $D_f$ containing $f$ such that $D_f\subseteq B\cup f$. The proof of (ii) may be obtained by switching $A$ and $B$. $\qed$
\bigskip

\noindent {\bf Lemma 4.3.} {\it Let $N$ be a simple and cosimple matroid in $\mathcal{M}$ with an exact $k$-separation $(A, B)$. Suppose $M\in \mathcal{M}$ is a simple and cosimple matroid such that $M\backslash e/f=N$. \begin {enumerate}
\item [(i)] If $\lambda_{M/f}(A)=k-1$  and $\lambda_{M\backslash e}(A)=k-1$, then $\lambda_M(A)=k-1$. 
\item [(ii)] If $\lambda_{M/f}(A\cup e)=k-1$  and $\lambda_{M\backslash e}(A\cup f)=k-1$, then $\lambda_M(A\cup\{e,f\})=k-1$. 
\end {enumerate}}

\bigskip

\noindent {\bf Proof.} Observe that, $M/f$ is a single-element extension of $N$ by element $e$ and $M\backslash e$ is a single-element coextension of $N$ by element $f$ (see Figure 5). 
Since $\lambda_{M/f}(A)=k-1$, Lemma 4.1(i) implies that     $M/f$ has a circuit $C_e$ such that $e\in C_e\subseteq B\cup e$. 
Since $\lambda_{M\backslash e}(A)=k-1$, Lemma 4.2(i) implies that  $M\backslash e$ has a cocircuit $D_f$ such that $f\in D_f\subseteq B\cup f$.  
It follows that, $$r_M(A)=r_{M\backslash e}(A)=r_N(A)$$  and
$$r_M(B\cup \{e, f\})=r_{M/f}(B\cup e)+1=r_N(B)+1.$$ 
Therefore,
\begin{eqnarray*}
 \lambda _M\left(A\right) & = & r_M  \left(A\right)  +r_M\left(B\cup \{e, f\}\right) - r \left(M\right)    \\
                  & = & r_{M\backslash e}\left(A\right)+r_{M/f}\left(B\cup e\right)+1- r\left(N\right)-1\\
                  & = & r_N\left(A\right)+r_N\left(B\right) - r\left(N\right)\\
                  & = & \lambda_N(A) \\  
                  & = & k-1    
\end{eqnarray*} 

The proof of (ii) may be obtained by switching $A$ and $B$. $\qed$
         
\bigskip

\begin{figure}[h]
\centering
\epsfxsize 2in \epsfbox{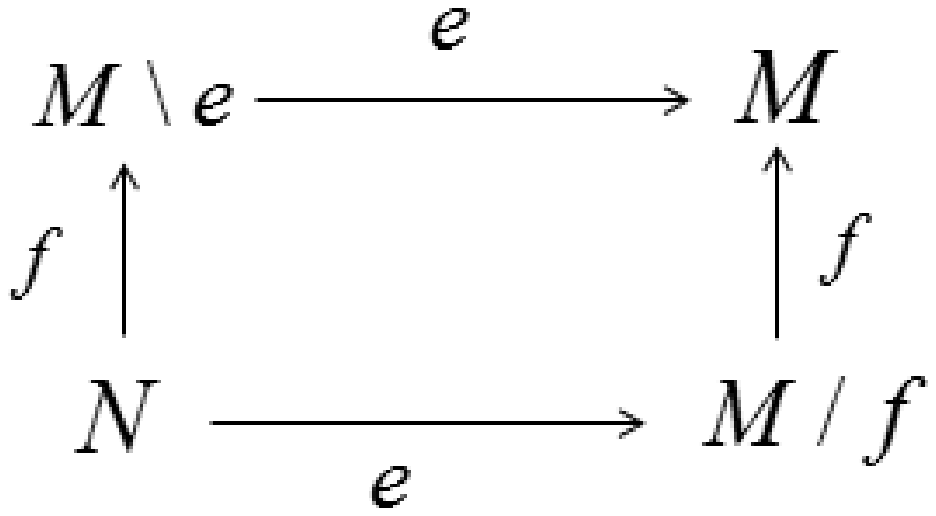}
\caption{$M\backslash e / f=N$ }
\end{figure}
\bigskip
 
We will also use standard properties of circuits and cocircuits. Suppose $M$ is a simple and cosimple matroid with the property that $M/f=N$, where $N$ has no loops nor coloops. Let $C$ be a circuit of $M$. If $f\in C$, then $C-f$ is a circuit in $M/f$. If $f\not\in C$, then $C$ is a circuit or union of circuits  in $M/f$. Looking at things the other way around, let $C'$ be a circuit in $N$. Then either $C'$ remains a circuit in $M$ or $C'\cup f$ is a circuit in $M$.  

Likewise for cocircuits. Suppose $M$ is a simple and cosimple matroid with the property that $M\backslash e=N$, where $N$ has no loops nor coloops. Let $D$ be a cocircuit of $M$. If $e\in D$, then $D-e$ is a cocircuit in $M\backslash e$. If $e\not\in D$, then $D$  is a cocircuit or union of cocircuits in $M\backslash e$. Let $D'$ be a circuit in $N$. Then either $D'$ remains a cocircuit in $M$ or $D'\cup e$ is a cocircuit in $M$.  

\bigskip

\noindent {\bf Proof of Theorem 1.1.} We  prove that Theorem 1.1 is equivalent to Theorem 2.2. Lemmas 3.1 and 3.2 show that Theorem 1.1(i) and (ii) are equivalent to Theorem 2.2(i) and (ii), respectively. We must show that Theorem 1.1(iii) is equivalent to Theorem 2.2(iii). First, assume Theorem 2.2(iii) holds. Let us consider four cases:

\medskip
\noindent {\it Case 1.} Suppose $\lambda _{M/f}(A)=k-1$ and $\lambda _{M\backslash e}(A)=k-1$. Then Lemma 4.3(i) implies that $\lambda _M(A)=k-1$. By Lemma 4.1(i) and Lemma 4.2(i) this happens if and only if $M$ has a circuit $C_e$ such that $e\in C_e\subseteq B\cup e$ and a cocircuit $D_f$ such that $f\in D_f\subseteq B\cup f$. The conditions in Theorem 2.2(iii) are satisfied.

\medskip
\noindent {\it Case 2.} Suppose $\lambda _{M/f}(A)=k-1$ and $\lambda _{M\backslash e}(A\cup f)=k-1$. Then, by hypothesis $\{e, f, g\}$ is a triad or triangle with $g\in A$ or $\lambda _M(A\cup f)=k-1$. In the first case there is nothing to show since the conditions in Theorem 2.2(iii) are satisfied. Therefore, suppose $\lambda _M(A\cup f)=k-1$. Then since $\lambda _{M\backslash e}(A\cup f)=k-1$ and $\lambda _M(A\cup f)=k-1$ and $M$ is a simple single-element extension of $M\backslash e$, Lemma 4.1(i) implies that $M$ has a circuit $C_e$ such that $e\in C_e\subseteq B\cup e$. So Theorem 2.2(iii)(b) is satisfied. 

Since $\lambda _{M/f}(A)=k-1$ and $\lambda _M(A\cup f)=k-1$ and $M$ is a cosimple single-element coextension of $M/f$, Lemma 4.2 (ii) implies that $M$ has a cocircuit $D_f$ such that $f\in D_f\subseteq A\cup f$. Now, suppose, if possible $M$ has a circuit $R$ such that $\{e, f\}\subseteq R\subseteq B\cup \{e, f\}$. Then $R\cap D_f=\{f\}$; a contradiction. Therefore, Theorem 2.2(iii)(a) is vacuously satisfied.

\medskip
\noindent {\it Case 3.} Suppose $\lambda _{M/f}(A\cup e)=k-1$ and $\lambda _{M\backslash e}(A)=k-1$. Then by hypothesis $\{e, f, g\}$ is a triad or triangle with $g\in A$ or $\lambda _M(A\cup e)=k-1$. In the first case there is nothing to show as in Case 2. Therefore, suppose $\lambda _M(A\cup e)=k-1$. Then since $\lambda _{M/f}(A\cup e)=k-1$ and $\lambda _M(A\cup e)=k-1$ and $M$ is a cosimple single-element coextension of $M/f$, Lemma 4.2(i) implies that $M$ has a cocircuit $D_f$ such that $f\in D_f\subseteq B\cup f$. So Theorem 2.2(iii)(a) is satisfied. 

Since $\lambda _{M\backslash e}(A)=k-1$ and $\lambda _M(A\cup e)=k-1$ and $M$ is a simple single-element extension of $M\backslash e$, Lemma 4.1(ii) implies that $M$ has a circuit $C_e$ such that $e\in C_e\subseteq A\cup e$. Now, suppose, if possible $M$ has a cocircuit $R$ such that $\{e, f\}\subseteq R\subseteq B\cup \{e, f\}$. Then, $R\cap C_e=\{e\}$; a contradiction. Therefore, Theorem 2.2(iii)(b) is vacuously satisfied.

\medskip
\noindent {\it Case 4.} Finally, suppose $\lambda _{M/f}(A\cup e)=k-1$ and $\lambda _{M\backslash e}(A\cup f)=k-1$. Then, by hypothesis $\{e, f, g\}$ is a triad or triangle in $M$ with $g\in A$, which satisfies the conditions in Theorem 2.2(iii).
\medskip

Conversely, assume Theorem 2.2(iii) holds.  Again, let us consider four cases:
\medskip

\noindent {\it Case 1.} Suppose $M$ has a cocircuit $D_f$ such that $f\in D_f\subseteq B\cup \{e, f\}$ and a circuit $C_e$ such that $e\in C_e\subseteq B\cup \{e, f\}$. Then $M\backslash e$ has a cocircuit $D_f'$ such that $f\in D_f'\subseteq B\cup f$. By Lemma 4.1(ii) this occurs if and only if $\lambda _{M\backslash e}(A)=k-1$. 
Similarly, $M/f$ has a circuit $C_e'$ such that $e\in C_e'\subseteq B\cup e$. By Lemma 4.2(ii) this occurs if and only if $\lambda _{M/f}(A)=k-1$.  So Theorem 1.1(iii)(a) is satisfied.

\medskip
\noindent {\it Case 2.} Suppose $M$ has a cocircuit $D_f$ such that $f\in D_f \subseteq B\cup \{e, f\}$, but no circuit $C_e$ such that $e\in C_e\subseteq B\cup \{e, f\}$ and $\{e, f, g\}$ is not a triangle with $g\in A$. Then, by an argument similar to Case 1, $\lambda _{M\backslash e}(A)=k-1$ and $\lambda _{M/f}(A)\neq k-1$. We may assume that 
$\lambda _{M/f}(A\cup e) = k-1$.

By hypothesis, there is no cocircuit $R$ such that $\{e, f\}\subseteq R\subseteq B\cup \{e, f\}$. So the cocircuit $D_f$ is such that $f\in D_f\subseteq B\cup f$. 
It follows from Lemma 4.1(ii) that $\lambda _{M/f}(A\cup e)=k-1$ if and only if $M/f$ has a circuit $C_e$ such that $e\in C_e\subseteq A\cup e$.
Suppose $C_e\cup f$ is a circuit in $M$. Then $C_e\cap D_f=\{f\}$; a contradiction. So $C_e$ stays a circuit in $M$. It follows that $\lambda _M(A\cup e)=k-1$, and the condition in Theorem 1.1(iii)(c) is satisfied.

\medskip
\noindent {\it Case 3.} Suppose $M$ has no cocircuit $D_f$ such that $f\in D_f \subseteq B\cup \{e, f\}$, and $\{e, f, g\}$ is not a triad with $g\in A$, but $M$ has a circuit $C_e$ such that $e\in C_e\subseteq B\cup \{e, f\}$. Then, by an argument similar to Case 1, $\lambda _{M\backslash e}(A)\neq k-1$ and $\lambda _{M/f}(A)= k-1$. We may assume that 
$\lambda _{M\backslash e}(A\cup f)= k-1$.

By hypothesis, there is no circuit $R$ such that $\{e, f\}\subseteq R\subseteq B\cup \{e, f\}$. So the circuit $C_e$ is such that $e\in C_e\subseteq B\cup e$. 
It follows from Lemma 4.2(ii) that $\lambda _{M\backslash e}(A\cup f)=k-1$ if and only if $M\backslash e$ has a cocircuit $D_f$ such that $f\in D_f\subseteq A\cup f$.
Suppose $D_f\cup e$ is a cocircuit in $M$. Then $D_f\cap C_e=\{e\}$; a contradiction. So $D_f$ stays a cocircuit in $M$. It follows that $\lambda _M(A\cup f)=k-1$ and Theorem 1.1(iii)(b) is satisfied.
 
\medskip

\noindent {\it Case 4.} Suppose $M$ has no circuit $C_e$ such that $e\in C_e\subseteq B\cup \{e, f\}$ and no cocircuit $D_f$ such that $f\in D_f \subseteq B\cup \{e, f\}$. Then by an argument similar to Case 1, $\lambda _{M\backslash e}(A)\neq k-1$ and $\lambda _{M/f}(A)\neq  k-1$. We may assume that $\lambda _{M\backslash e}(A\cup f)= k-1$ and $\lambda _{M/f}(A\cup e)= k-1$. Lemma 4.3(ii) implies that $\lambda _M(A\cup \{e, f\})=k-1$. However, in this situation, the exact 3-separation is not necessarily maintained in 3-connected matroids containing $M$, as shown by the example at the end of the next section. So, in this case the only possibility is for $\{e, f, g\}$ to be a triangle or triad with $g\in A$. Hence the condition in Theorem 1.1(iii)(d) is satisfied.
 $\qed$

\bigskip

\noindent {\bf Proof of Corollary 1.2. } Since $A$ is a $4$-element circuit and cocircuit in $N$ and its simple single-element extensions and cosimple single-element coextensions, $\lambda (A)=2$ in all these matroids. Hence, the conditions in Theorem 1.1(i) and (ii) are immediately satisfied. Suppose  $M\backslash e/f=N$ and $M\backslash e$ is simple and $M/f$ is cosimple, then the result follows from Theorem 1.1(iii)(a). It remains to consider the case when either $M\backslash e$ is not cosimple or $M/f$ is not simple.  

Suppose $M\backslash e$ is not cosimple. Then $\{e, f, g\}$ is a triad in $M$ and $f$ is in series to an element $g$ in $M\backslash e$. If $g\in A$,   then Theorem 1.1 (iii)(b) holds. If $g\in B$,  there is a cocircuit  $\{f, g\}\subset B\cup f$  containing $f$. By Lemma 4.2(i), $\lambda _{M/f}(A)=2$ and Theorem 1.1(iii)(a) is satisfied.

Suppose $M/f$ is not simple. Then $\{e, f, g\}$ is a triangle in $M$ and $f$ is parallel to an element $g$ in $M/f$. If $g\in A$,   then Theorem 1.1(iii)(c) holds. If $g\in B$, then there is a circuit  $\{f, g\}\subset B\cup f$  containing $f$. By Lemma 4.1(i), $\lambda _{M/f}(A)=2$ and Theorem 1.1(iii)(a) is satisfied.  $\qed$
\bigskip


\section {$R_{12}$ is a $3$-decomposer for regular matroids}

Recall that,  $R_{12}$ has a non-minimal exact 3-separation $(A, B)$, where $A=\{3, 4, 7, 8, 11, 12\}$ and $B=\{1, 2, 5, 6, 9, 10\}$. Both $A$ and $B$ are the union of circuits and the union of cocircuits (this is just a coincidence). Using Theorem 1.1, the biggest reduction in computation occurs when computing $M$, where $M\backslash e/f=N$ because if $\lambda _{M/f}(A)=2$ and $\lambda _{M\backslash e}(A)=2$, there is no need to check $M$.

\bigskip
\noindent {\bf Theorem 4.1.} {\it $R_{12}$ is a $3$-decomposer for regular matroids.}

\bigskip

\noindent {\bf Proof.}  $R_{12}$ has two non-isomorphic single-element extensions that are regular, namely, $P_{13}$ and $Q_{13}$, where $P_{13}$ is formed by adding any one of the columns $\alpha=[0 0 0 0 1 1]$, $\beta=[1 1 0 0 0 0]$, or $\gamma=[1 1 0 0 1 1]$ and $Q_{13}$ is formed by adding column $\delta=[0 0 1 1 0 0]$. We can check that $\lambda (A)=2$ for the first three columns, but for the last column $d=[0 0 1 1 0 0]$, $\lambda (A\cup e)=2$ 

Since $R_{12}$ is self-dual (and this particular representation is symmetric) $R_{12}$ has two non-isomorphic single-element coextensions that are regular, namely $P_{13}^*$ and $Q_{13}^*$, where $P_{13}^*$ is formed by adding any one of the rows $a=[0 0 0 0 1 1]$, $b=[1 1 0 0 0 0]$, or $c=[1 1 0 0 1 1]$ and $Q_{13}^*$ is formed by adding row $d=[0 0 1 1 0 0]$. We can check that $\lambda (A)=2$ for all four rows.  Therefore, the conditions in Theorem 1.1(i) and (ii) are met.

Next, suppose $M$ is a cosimple single-element coextension of a Type (i) matroid such that $M\backslash e/f = N$. By Theorem 1.1(iii) the only $M$ whose connectivity condition we must check is $M$ where $M/f$ is formed by the column $[0 0 1 1 0 0]$. In other words, we must check the cosimple coextensions of $Q_{13}$ formed by adding rows  $a$, $b$, $c$, or $d$ with a zero or a one at the end and the rows corresponding to entries in $B$ with the last entry reversed, namely, $1'=[1 1 1 0 0 1]$, $2'=[1 1 0 1 0 0 1]$, $5'=[0 0 1 0 1 1 0]$, $6'=[0 0 0 1 1 1 1]$, $10'=[0 0 1 0 0 0 1]$, $11'=[0 0 0 1 0 0 1]$ because for these rows $\lambda _{M\backslash e}(A)=2$. Table 1a shows that in all cases where the matroid is regular, $\lambda _M(A\cup e)=2$. So Theorem 1.1(iii) is satisfied. 

Lastly, suppose $M$ is a simple single-element extension of a Type (ii) matroid such that $M\backslash e/f = N$. By Theorem 1.1(iii) we must check every coextension of $R_{12}$, but the only columns we have to add to the coextension are  $[0 0 1 1 0 0 0]$ and $[0 0 1 1 0 0 1]$. Table 1b shows that in all cases where the matroid is regular, $\lambda _M(A\cup e)=2$. So Theorem 1.1(iii) is satisfied.  $\qed$
\bigskip

\small
 \begin{center}
\begin{tabular}{|c|c|c|c|}
\hline

\bf   &\bf{Coext. Rows} & \bf{$F_7$ or $F_7^*$-minor} & \bf{3-separation}  or triad \\  \hline \hline

$Q_{13}$ with column $[0 0 1 1 0 0]$  & $[0 0 0 0 1 1 0]$  & No &   $\lambda\{3, 4, 8, 9, 12, 13, {\bf 14}\}=2$  \\ \hline

& $[0 0 0 0 1 1 1]$  & YES &     \\ \hline
                                      
& $[1 1 0 0 0 0 0]$  & No &  $\lambda\{3, 4, 8, 9, 12, 13, {\bf 14}\}=2$  \\ \hline

& $[1 1 0 0 0 0 1]$  & YES &   \\ \hline
  
& $[1 1 0 0 1 1 0]$  & No &  $\lambda\{3, 4, 8, 9, 12, 13, {\bf 14}\}=2$  \\ \hline

& $[1 1 0 0 1 1 1]$  & YES &    \\ \hline

& $[0 0 1 1 0 0 0]$  & No &  $\lambda\{3, 4, 8, 9, 12, 13, {\bf 14}\}=2$  \\ \hline

& $[0 0 1 1 0 0 1]$  & YES & \\ \hline  

& $[1 1 1 0 0 0 1]$  & YES & \\ \hline  
& $[1 1 0 1 0 0 1]$  & YES & \\ \hline 
& $[0 0 1 0 1 1 0]$  & YES & \\ \hline  
& $[0 0 0 1 1 1 1]$  & YES & \\ \hline  
& $[0 0 1 0 0 0 1]$  & YES & \\ \hline 
& $[0 0 0 1 0 0 1]$  & YES &  \\ \hline\hline 

\end{tabular}
\end{center}
\normalsize
 \begin{center} Table 1a: Regular cosimple single-element coextensions of $Q_{13}$  \end{center} 

\bigskip

\small
 \begin{center}
\begin{tabular}{|c|c|c|c|}\hline

\bf   &\bf{Ext. Columns} & \bf{$F_7$ or $F_7^*$-minor} & \bf{3-separation}   \\  \hline \hline

$P_{13}^*$ with row $[0 0 0 0 0 1 1]$  & $[0 0 1 1 0 0 0]$  &  No  &   $\lambda\{3, 4, 8, 9, 12, 13, {\bf 14}\}=2$  \\ \hline
 & $[0 0 1 1 0 0 0 1]$  & YES  &     \\ \hline \hline 

$P_{13}^*$ with row $[1 1 0 0 0 0]$  & $[0 0 1 1 0 0 0]$   &  No &   $\lambda\{3, 4, 8, 9, 12, 13, {\bf 14}\}=2$  \\ \hline
& $[0 0 1 1 0 0 1]$  & YES &     \\ \hline \hline 

$P_{13}^*$ with row $[1 1 0 0 0 1 1]$  & $[0 0 1 1 0 0 0]$  &  No &   $\lambda\{3, 4, 8, 9, 12, 13, {\bf 14}\}=2$  \\ \hline
& $[0 0 1 1 0 0 1]$  & YES &     \\ \hline \hline 

$Q_{13}^*$ with row $[0 0 1 1 0 0]$  & $[0 0 1 1 0 0 0]$  & No  &   $\lambda\{3, 4, 8, 9, 12, 13, {\bf 14}\}=2$  \\ \hline
& $[0 0 1 1 0 0 1]$   & YES &     \\ \hline \hline 

\end{tabular}
\end{center}
\normalsize
 \begin{center} Table 1b: Regular simple single-element extensions of $P_{13}^*$ and $Q_{13}^*$  \end{center} 

 \bigskip

Observe that $R_{12}$ is a very simple example of a 3-decomposer. It has few regular extensions and coextensions and not only is it self-dual, but the representation we are using is symmetric. Additionally, both $A$ and $B$ are union of circuits and union of cocircuits and as such interchangable. All this makes it an easy computational example. However, the algorithmic approach developed in this paper makes it possible to find more complicated 3-decomposers when there are many extensions and coextensions in the excluded minor class. In fact, it reduces the entire process of finding a 3-decomposer to the click of a button. This method will allow us to find 3-decomposers when it is no longer feasible to check the computations by hand or even display the tables as shown here.

Finally, we give an example to show that there is no symmetry between $A$ and $B$ that would allow us to conclude in step Theorem 1.1(iii) that $\lambda _M(A\cup \{e, f\})=k-1$ is acceptable, except when $\{e, f, g\}$ is a triad with $g\in A$. Consider the 10-element rank-5 matroid $M[X]$ shown below. It has a non-minimal exact 3-separation $(A, B)$ where $A=\{1, 2, 5, 6, 7, 10\}$ and $\lambda_X (A)=2$.

\[
X=\left[ 
\begin{array}{ccccc|ccccc}
&&&&&    0&1&1&1&1 \\
&&&&&    1&0&1&1&1 \\
&&I_5&&& 1&1&0&1&0 \\
&&&&&    1&1&1&1&0 \\
&&&&&    0&1&0&0&1
\end{array} 
\right] 
\]

If we extend $X$ by column $[1 1 0 0 0]$, we get a matroid $Y$ with the property that $\lambda _Y(A)=2$. Now, let us coextend $Y$ by $[1 1 0 0 1 1]$ to get the matroid $Z$. In $Z$, $\lambda _Z(A)\neq 2$, $\lambda _Z(A\cup \{6\})\neq 2$ and there is no triad $\{e, f, g\}$  with $g\in A$, but $\lambda _Z(A\cup \{6, 12\})=2$. 
Next, let us coextend $Y$ by row $[1 1 1 0 0 1]$ to get the matroid $Z'$. In $Z'$, $\lambda _{Z'}(A)=2$. 
However, if we coextend $Y$ by both rows  $[1 1 0 0 1 1]$ and $[1 1 1 0 0 1]$ we get the matroid $Q_{13}$ shown below that is internally 4-connected. 

\[
Q_{13}=\left[ 
\begin{array}{ccccc|cccccc}
&&&&&    0&1&1&1&1&1 \\
&&&&&    1&0&1&1&1&1 \\
&&I_7&&& 1&1&0&1&0&0 \\
&&&&&    1&1&1&1&0&0 \\
&&&&&    0&1&0&0&1&0 \\
&&&&&    1&1&0&0&1&1 \\
&&&&&    1&1&1&0&0&1 \\
\end{array} 
\right] 
\]

In terms of Theorem 2.2, the matroid $Z$ has a circuit $R$ such that $\{12, 6\}\subseteq R \subseteq B\cup \{12, 6\}$, but no cocircuit $D$ such that $6\in D\subset B\cup \{12, 6\}$ and there is no triad of the form $\{12, 6, g\}$ with $g\in A$. As noted earlier, this example also serves to show that a 2-element check for the presence of non-minimal exact 3-separations will not do for binary matroids. The conditions in Seymour's Decomposition Theorem and consequently in Theorem 1.1 are best possible.

\bigskip
\noindent {\bf Acknowledgements:} The author thanks James Oxley for being a mentor for over 20 years and Manoel Lemos for assistance with the proof of Theorem 2.2. 
\bigskip

\noindent {\bf References}
\bigskip

\begin{enumerate}

\item  D. Mayhew, G. Royle, and G. Whittle (2011)  The internally 4-connected binary matroids with no $M(K_{3,3})$-minor. {\it Memoirs of the American Mathematical Society} {\bf 981}, American Mathematical Society, Providence, Rhode Island.

\item J. G. Oxley (1992) {\it Matroid Theory}, Second Edition (2011), Oxford University Press, New York.
                                                                                                          
\item P. D. Seymour (1980) Decomposition of regular matroids, {\it J. Combin. Theory  Ser. B} {\bf 28}, (1980) 305-359.
\end{enumerate}
\bigskip

\end {document}